\documentclass[a4paper, 11pt]{article} 
\usepackage{amsfonts,amsmath,amssymb,amscd,color,stmaryrd,latexsym,shadow}
\usepackage[all]{xy} 
\usepackage[T1]{fontenc}
\usepackage[latin1]{inputenc}

\def\End{\mathop{\rm End}\nolimits}

\def\deg{\mathop{\rm deg}\nolimits}

\def\diff{\mbox{\sl Diff}}

\def\ric{\mathop{\rm Ric}\nolimits}

\def\note#1{\marginpar{\raggedright\if@twoside\ifodd\c@page\raggedleft\fi\fi\sf\scriptsize RMK: #1}}

\newcommand{\mf}{\mathfrak}
\newcommand{\mb}{\mathbf}
\newcommand{\mc}{\mathcal}

\newcommand{\R}{\mathbb{R}}
\newcommand{\Z}{\mathbb{Z}}
\newcommand{\C}{\mathbb{C}}

\newcommand{\E}{TM\oplus T^*\!M}

\newenvironment{remark}{\begin{trivlist}\item[]{\bf Remark:}\setlength{\parindent}{0pt}}{\end{trivlist}}
\newenvironment{example}{\begin{trivlist}\item[]{\bf Example:}\setlength{\parindent}{0pt}}{\end{trivlist}}
\newenvironment{examples}{\begin{trivlist}\item[]{\bf Examples:}\setlength{\parindent}{0pt}}{\end{trivlist}}

\newtheorem{theorem}{Theorem}[section]
\newtheorem{definition}[theorem]{Definition}
\newtheorem{proposition}[theorem]{Proposition}

\newtheorem{corollary}[theorem]{Corollary}

\begin{document}
\thispagestyle{empty}
\vspace{12pt}

\begin{center}{\LARGE{\bf Calabi--Yau manifolds with $B$--fields}}

\vspace{35pt}

{\bf Frederik Witt\footnote{The present article is an extended version of a talk given at the Universit\`a di Torino and the Erwin--Schr\"odinger Institute (ESI), Wien. I wish to thank A.~Fino and S.~Console for the very enjoyable time spent at Torino. I also like to thank the organisers of the workshop ``Poisson Sigma Models, Lie Algebroids, Deformations, and Higher Analogues'' and the Erwin Schr\"odinger Institute for hospitality whilst writing this paper. Further thanks are due to C.~Jeschek and F.~Gmeiner from whom I learned a great deal of physics, and to D.~Ploog for proofreading. Finally, I am indebted to the supervisor of my doctoral thesis, N.~Hitchin, who introduced me to this subject. The author was supported by the DFG as a member of the SFB 647.}
}

\vspace{15pt}

{\it NWF I - Mathematik\\
Universit\"at Regensburg\\
D--93040 Regensburg, F.R.G.\\
e-mail: \texttt{Frederik.Witt@mathematik.uni-regensburg.de}}\\[1mm]

\vspace{40pt}

{\bf ABSTRACT}
\end{center}
In recent work N.~Hitchin introduced the concept of ``generalised geometry''. The key feature of generalised structures is that that they can be acted on by both diffeomorphisms and $2$--forms, the so--called $B$--fields. In this lecture, we give a basic introduction and explain some of the fundamental ideas. Further, we discuss some examples of generalised geometries starting from the usual notion of a Calabi--Yau manifold, as well as applications to string theory.

\bigskip

%
%
%
%
%
\section{Introduction}
This text tries to motivate and to explain a new type of geometry, the so--called {\em generalised geometry}, which was developed in Hitchin's foundational article~\cite{hi03} and taken further by his students~\cite{ca05},~\cite{gu04},~\cite{wi05}.

\bigskip

In differential geometry, one is familiar with deforming a geometric structure by a diffeomorphism. For instance, if $(M,g)$ is a Riemannian manifold, then so is $(M,f^*g)$ for any $f\in\diff(M)$. Enlarging the notion of a Riemannian structure to a generalised Riemannian structure consists in allowing deformations by a $2$--form, or, following standard physicists' jargon, by a {\em $B$--field}. This concept makes also sense for other geometries. The most special class we will consider in this note is that of a generalised $SU(m)$--structure which comprises classical Calabi--Yau manifolds and their $B$--field deformations. Non--metric geometries such as symplectic structures can also be generalised, and we will encounter examples hereof in the sequel.

\bigskip

The motivation for looking at $B$--fields originated from physics. Not surprisingly then, the most striking applications of generalised geometry so far arose in string theory, where generalised structures accounted for the geometry imposed by supersymmetry. Despite the resulting neat formulation and new insight gained, one might question the use of general geometry in mathematics. The answer is, as often, {\em generality}: Generalised geometry unifies different concepts and combines structures as different such as special complex and symplectic manifolds within one geometric framework. Further, as already emphasised in~\cite{gu04}, generalised geometry is naturally described in terms of Courant algebroids which links these to other areas such as foliated or Poisson geometry.

\bigskip

The plan of this note is as follows. We start by summing up some key results of the structures we will generalise in the sequel such as complex, K\"ahler and Calabi--Yau manifolds. We then explain the basic idea behind the concept of supersymmetry and how this involves the special geometric structures we discussed at the beginning. This also leads to investigate $B$--field actions which are considered in Section~\ref{btrafo}. Section~\ref{gengeo} contains the core part where we describe explicitly the generalisations of the geometric structures discussed in Section~\ref{classgeom}. More details and proofs can be mainly found in~\cite{gu04},~\cite{hi03} and~\cite{wi05} upon which this text is based. It is not intended to give a complete overview, and several interesting topics had to be glossed over among which there are, for instance, moduli spaces of generalised $K3$--surfaces and their relation to $\mc{N}=(2,2)$ super conformal field theory~\cite{hu05}, topological strings 
 and the generalised $A$-- and $B$--model~\cite{cgj06},~\cite{ka04}~\cite{pe07}, generalised versions of the $\partial\overline{\partial}$--lemma~\cite{ca05} or generalised calibrations~\cite{gmwi06},~\cite{gmwi07} to mention only a few articles and researchers who contributed to the field. At any rate, I hope to convince the reader of the usefulness of the concepts to be presented below.
%
%
%
%
%
\section{Classical geometries}
\label{classgeom}
Let $M^n$ be a smooth, orientable manifold of dimension $n$. If $n$ is even, we write $n=2m$. By a {\em classical geometry} or {\em structure}, we understand a collection of tensors over the tangent bundle $TM$ enjoying special algebraic properties. Further, these tensors are possibly subject to integrability conditions. An alternative viewpoint is provided by $G$--structure theory (cf. the classical text~\cite{kono63} or~\cite{sa89}), where geometric structures are described in terms of reductions of the principal $GL(n)$--frame bundle to a principal $G$--fibre bundle. One also says that $TM$ is {\em associated} with a $G$--structure and refers to $G$ as the structure group of $TM$. Locally, after choosing some trivialisation $U\times\R^n\cong TM_{|U}$ of the tangent bundle, we can think of a $G$--structure as a smooth family of Lie group representations $r_x:G\to GL(T_xM)$ parametrised by points $x\in U$ which is well behaved under transitions of two different trivialisatio
 ns. 
\paragraph{Riemannian structures.}
A {\em Riemannian manifold} $(M,g)$ is specified by a symmetric $(2,0)$--tensor field $g$ whose special algebraic property is to be {\em positive definite}, i.e.\ $g(X,X)>0$ for $X\in TM$ other than $0$. This endows every tangent space $T_xM$ with the structure of an oriented Euclidean vector space. The symmetry group, that is, the stabiliser of this datum under the action of $GL(T_xM)$, is $SO(T_xM,g_x)$. Locally, we thus obtain a family of representations $r_x:SO(n)\stackrel{\sim}{\to} SO(T_xM,g_x)\subset GL(T_xM)$, and since $g$ is a globally defined object, these families come from a principal $SO(n)$--fibre bundle. Therefore, a Riemannian structure amounts to saying that $TM$ is associated with $SO(n)$. 
\paragraph{Symplectic structures.}
An {\em almost symplectic manifold} $(M^{2m},\omega)$ is specified by a skew--symmetric $2$--form $\omega$ whose special algebraic property is to be {\em non--degenerate}, that is, $\omega^m\not=0$. The symmetry group is $Sp(m,\R)$, so locally we obtain families of representations $r_x:Sp(m,\R)\stackrel{\sim}{\to} Sp(T_xM,\omega)\subset GL(T_xM)$. As for Riemannian structures, we see that an almost symplectic structure is nothing else than an $Sp(m,\R)$--structure. This structure is called {\em integrable} or simply {\em symplectic} if $d\omega=0$, i.e.\ $\omega$ is closed. Darboux's theorem asserts that this is equivalent to finding local coordinates $p^1,\ldots,p^m,q^1,\ldots,q^m$ on $U\subset M$ such that $\omega_{|U}=\sum dp^j\wedge dq^j$. 
\paragraph{Complex structures.}
An {\em almost complex manifold} $(M^{2m},J)$ is specified by a $(1,1)$--tensor field $J$, that is, a bundle morphism $J:TM\to TM$. Its special algebraic property is to square to minus the identity, $J^2=-Id$. This endows every tangent space with the structure of a complex vector space, and locally, we have representations $r_x:GL(m,\C)\stackrel{\sim}{\to} GL(T_xM,J_x)\subset GL(T_xM)$, giving rise to a $GL(m,\C)$--structure. Equivalently, we have a decomposition of the complexification $TM\otimes\C=T^{1,0}\oplus T^{0,1}$ into complex vector bundles which are conjugate to each other, $\overline{T^{1,0}}=T^{0,1}$, namely the $\pm i$--eigenbundles of $J$. One speaks of an {\em integrable almost complex structure} or simply a {\em complex structure}, if the associated {\em Nijenhuis tensor} 
$$
N^J(X_0,X_1)=[X_0,X_1]-[JX_0,JX_1]+J([JX_0,X_1]+[X_0,JX_1])
$$ 
vanishes. Equivalently, $T^{1,0}$ (or $T^{0,1}$) is preserved under the extension of the Lie bracket to $TM\otimes\C$. The dual version is this: Let $T^*\!M\otimes\C=T^{1,0*}\oplus T^{0,1*}$ be the decomposition induced by the dual map $J^*:T^*\!M\to T^*\!M$. Put $\Lambda^{p,q}=\Lambda^pT^{1,0*}\otimes\Lambda^qT^{0,1*}$, whence $\Lambda^kT^*\!M\otimes\C=\oplus_{p+q=k}\Lambda^{p,q}$. We define differential operators 
$$
\partial:\Lambda^{p,q}\to\Lambda^{p+1,q},\quad\overline{\partial}:\Lambda^{p,q}\to\Lambda^{p,q+1}
$$
by compounding the exterior differential with the corresponding projection operator. As this can be computed in terms of vector field derivations and the Lie bracket, it follows that a complex structure is integrable if and only if $d=\partial+\overline{\partial}$. Integrability is also equivalent to the existence of complex coordinates $z^1=x^1+iy^1,\ldots, z^m=x^m+iy^m$, that is, $J(\partial_{x^j})=\partial_{y^j}$ and $J(\partial_{y^j})=-\partial_{x^j}$. The bundles $T^{1,0}$ and $T^{0,1}$ are spanned by $\partial_{z^j}=\partial_{x^j}-i\partial_{y^j}$ and $\partial_{\overline{z}^j}=\partial_{x^j}+i\partial_{y^j}$ respectively, while the differentials $dz^j=dx^j+idy^j$ and $d\overline{z}^j=dx^j-idy^j$ span $T^{1,0*}$ and $T^{0,1*}$. 

\bigskip

An important subclass of complex structures is provided by {\em special complex manifolds} $(M^{2m},\Omega)$, where $\Omega\in\Omega^m(M)\otimes\C$ satisfies the following properties: (1) $\Omega\wedge\overline{\Omega}\not=0$, (2) locally $\Omega=\theta^1\wedge\ldots\wedge\theta^m$ for complex $1$--forms $\theta^j$, i.e.\ $\Omega$ is locally decomposable, and (3) $d\Omega=0$. By (1) and (2), we obtain a splitting $T^*\!M\otimes\C=T^{1,0*}\oplus T^{0,1*}$, where the $(1,0)$--forms $T^{1,0*}$ are locally spanned by the $\theta^j$. Hence $\Omega\wedge\theta=0$ for any $(1,0)$--form $\theta$, so that (3) implies $\Omega\wedge d\theta=0$. Consequently, $d\theta$ has no $(0,2)$--part, which implies integrability via the exterior derivative characterisation. In particular, $\overline{\partial}\Omega=0$, that is, $\Omega$ defines a {\em holomorphic} trivialisation of the {\em canonical line bundle} $\kappa=\Lambda^{m,0}=\Lambda^mT^{1,0*}$. Note that $T^{1,0}$ is isomorphic with the t
 angent space $TM$ as a complex vector bundle. Since $T^{1,0}_x$ comes equipped with the complex volume form $\Omega_x$ whose stabiliser under the action of $GL_{\C}(T^{1,0}_x)$ is $SL_{\C}(T^{1,0}_x,\Omega_x)$, the structure group is $SL(m,\C)$ (whence the name of special complex). 
\paragraph{K\"ahler structures.} 
An {\em almost K\"ahler manifold} $(M,J,g)$ is an almost complex manifold which carries a hermitian metric $g$, that is, $J$ acts as an isometry for $g$. The structure group is thus $U(m)$. The structure is {\em integrable} or simply {\em K\"ahler}, if $J$ is integrable and the {\em K\"ahler form} defined by 
\begin{equation}\label{symcommet}
\omega(X_0,X_1)=g(JX_0,X_1)
\end{equation} 
is symplectic, i.e $d\omega=0$. Note that conversely, if $(M,\omega)$ defines a symplectic structure on the complex manifold $(M,J)$ such that $\omega(JX_0,JX_1)=\omega(X_0,X_1)$ and
\begin{equation}\label{positivity}
\omega(X,JX)>0
\end{equation} 
for $X\not=0$, then~(\ref{symcommet}) defines a K\"ahler metric. Hence, a K\"ahler manifold is equally specified by the triple $(M,J,\omega)$.
\paragraph{Calabi--Yau structures.} 
An {\em almost Calabi--Yau manifold} $(M^{2m},\omega,\Omega)$ is made up out of an almost symplectic structure $\omega$ and an almost special complex structure $\Omega$ subject to the compatibility conditions (1) $\Omega\wedge\overline{\Omega}=(-1)^{m(m-1)/2}i^m\omega^m$ and (2) $\Omega\wedge\omega=0$. Further, if $J$ denotes the induced complex structure, then $\omega$ has to verify the positivity condition~(\ref{positivity}). The associated structure group is $SU(m)$. In particular, the inclusion $SU(m)\subset U(m)$ induces an almost K\"ahler metric $g$. The structure is {\em integrable} or {\em Calabi--Yau} for short, if $\omega$ and $\Omega$ are closed. For $M$ compact, the Calabi--Yau theorem asserts the following converse: If $M$ is K\"ahler with vanishing first Chern class and K\"ahler form $\omega$, there exists a (unique) metric $g'$ induced by an underlying Calabi--Yau structure with $[\omega']=[\omega]\in H^2(M,\R)$. 

\bigskip

A different characterisation of Calabi--Yau manifolds can be given in terms of spinor fields. As this concept will be of much use later on, we briefly outline the general setting; good references are~\cite{ba81} and~\cite{lami89}. Let $V$ be an oriented real vector bundle of rank $k$ with a bundle metric  $(\cdot\,,\cdot)$ of signature $(p,q)$. A {\em Clifford bundle} for $V$ is a complex vector bundle $\mb{S}(V)\to M$ with a hermitian metric $h$ and a $C^{\infty}(M)$--linear action $\mu:\Lambda^*V\times\mb{S}(V)\to\mb{S}(V)$, the so--called {\em Clifford multiplication}, such that
$$
\mu\big(X,\mu(Y,\rho)\big)+\mu\big(Y,\mu(X,\rho)\big)=-(X,Y)\rho,
$$ 
and
$$
h\big(\mu(X,\rho),\tau\big)+h\big(\rho,\mu(X,\tau)\big)=0,
$$
where $X$ and $Y$ are sections of $V$. This definition extends to the exterior algebra by
$$
\mu(X\wedge Y,\rho)=\mu\big(X,\mu(Y,\rho)\big)
$$
for $X,\,Y\in V$ orthogonal etc. If $\nabla$ is a metric connection on $V$, then a hermitian connection $\nabla^{\mb{S}}$ on $\mb{S}(V)$ is a {\em Clifford connection} if
$$
\nabla^S\mu(X,\rho)=\mu(\nabla X,\rho)+\mu(X,\nabla^S\rho).
$$
If ${\rm rk}\big(\mb{S}(V)\big)=2^{[{\rm rk}(V)/2]}$, then $\mb{S}(V)$ is called a {\em (Dirac) spinor bundle}. It is associated with with a {\em spin structure}, which locally gives rise to a family of representations $r_x:Spin(p,q)\to U\big(\mb{S}(V)_x,h_x\big)\subset GL\big(\mb{S}(V)_x\big)$. Here $Spin(p,q)$ is the so--called {\em spin group of signature $(p,q)$}, the double cover of $SO(p,q)$. Sections of $\mb{S}(V)$ are called {\em spinor fields}. If the rank of $V$ is even, then there exists a natural decomposition $\mb{S}(V)=\mb{S}(V)_+\oplus\mb{S}(V)_-$ into spinor bundles of {\em positive} and {\em negative chirality}, where the fibres $\mb{S}(V)_{\pm x}$ are irreducible representations of $Spin(p,q)$. Sections of these are referred to as {\em positive} and {\em negative} respectively, or as {\em chiral} in general. The chirality is reversed under Clifford multiplication with $\Lambda^{od}V$. 

\bigskip

There are general existence and classification results for spin structures. In the cases we shall consider here, we can construct the spinor bundle explicitly in a canonical way. For instance, if $V=TM$ is endowed with an almost Calabi--Yau metric $g$, we find $\mb{S}(TM)_+=\Lambda^{0,ev}$ and $\mb{S}(TM)_-=\Lambda^{0,od}$. The hermitian structure comes from the hermitian metric on $T^{0,1*}$ induced by $g$. Further, Clifford multiplication is given by $\mu(X,\Psi)=-X\llcorner\Psi+g(\overline{X},\cdot)\wedge\Psi$, where we identified $TM$ with $T^{0,1}$ via the almost complex structure $J$. We will also write $X\cdot\Psi$ for $\mu(X,\Psi)$. In particular, $\Psi=1$ and $\Phi=\overline{\Omega}$ are two globally defined spinor fields of constant norm. They are related by the so--called {\em charge conjugation operator} $\mc{A}$, a conjugate linear operator which on each fibre commutes with the action of $Spin(2m)$. Further, $\Psi$ and $\Phi$ are {\em pure}, i.e.\ $W_{\Psi}=\{Z\in TM\otimes\C\,|\,Z\cdot\Psi=0\}$ is a complex subbundle of maximal rank $m$ ($W_{\Psi}$ being isotropic, the rank is at most $m$), and similarly for $W_{\Phi}$. A pure spinor field is necessarily chiral. The Levi--Civita connection $\nabla^g$ on $TM$ induces a canonical Clifford connection on $\mb{S}(TM)$ which preserves the splitting $\mb{S}(TM)_{\pm}$ and which we denote, by an {\em abus de notation}, also by $\nabla^g$. If the metric comes actually from an (integrable) Calabi--Yau, then $\nabla^g\Psi=0$ and $\nabla^g\Phi=0$ (note that this is not the usual covariant derivative on $\Lambda^{0,*}$, so $\Psi=1$ is a priori not covariant constant). Conversely, assume that $TM$ carries a Riemannian metric $g$ for which a spinor bundle $\mb{S}(TM)$ exists. The stabiliser of any pure spinor $\Psi_x\in\mb{S}(TM)_{\pm\,x}$ is isomorphic to $SU(m)$ (cf.~\cite{lami89}). Hence, any pure spinor field of unit norm induces an almost Calabi--Yau structure which is integrable if and only
  if $\nabla^g\Psi=0$. Since~$\mc{A}$ commutes with $\nabla^g$, the pure spinor $\mc{A}(\Psi)$ is also parallel. To make contact between the spinor and the form point of view, we recall that there is a canonical bundle isomorphism $\mb{S}(TM)\otimes\mb{S}(TM)\cong\Lambda^*T^*\!M\otimes\C$ under which
\begin{equation}\label{CYspinor}
\Psi\otimes\Psi=\Omega,\quad\mc{A}(\Psi)\otimes\Psi=e^{i\omega}
\end{equation}
(cf.~\cite{wa89}). Further, this isomorphism commutes with $\nabla^g$, so that $\nabla^g\Psi=0$ and $\nabla^g\mc{A}(\Psi)=0$ are equivalent to $\nabla^g\omega=0$ and $\nabla^g\Omega=0$. In turn, this is equivalent (although this is not entirely obvious) to $d\Omega=0$ and $d\omega=0$.

\begin{remark}
All the ``almost''--structures considered here were defined in terms of tensors. Therefore, they transform naturally under vector bundle isomorphisms of $TM$ (covering a bijective smooth map $M\to M$), as these preserve any special algebraic property of a tensor field. Further, if such an isomorphism is induced by a diffeomorphism, it commutes with the natural differential operators coming from the underlying differential structure such as the Lie bracket or exterior differentiation. In this case, integrability conditions are preserved, too. For instance, if $(M,\omega)$ is an almost symplectic manifold, its $f$--{\em transform} $(M,f^*\omega)$ for some vector bundle isomorphism $f$ of $TM$ also defines an almost symplectic structure. If $f$ comes from a diffeomorphism, then in addition $df^*\omega=f^*d\omega$, so the $f$--transform is integrable, i.e.\ symplectic, if and only if the original structure is integrable. Extending the diffeomorphism group will eventually lead to 
 so--called {\em generalised geometries}. Before explaining this, we shall make a detour into physics first.
\end{remark}
%
%
%
%
%
\section{Supersymmetry}
\label{physicsI}
To motivate the subsequent development, we briefly look at supersymmetric $\sigma$--models and supergravity compactifications. For a more detailed account see for instance~\cite{fgr98},~\cite{fr99} or~\cite{hal03}.

\bigskip
 
The key notion we shall discuss here is the concept of {\em supersymmetry}. To see where this originates from, recall the basic setup of quantum mechanics: We are given a (complex) Hilbert space $\mc{H}$ (whose projectivisation corresponds to the ``physical'' states), and a set of self--adjoint linear operators, the {\em observables} (corresponding to classical ``measureable'' quantities such as position, momentum, energy etc). In particular, the observable corresponding to energy is the {\em Hamiltonian} of the system. It turns out that $\mc{H}$ comes equipped with a natural $\Z_2$--grading which reflects the existence of two types of particles. Namely, there are {\em bosons} (particles that transmit forces such as photons), and {\em fermions} (particles that make up matter such as electrons). More generally, fields are either bosonic and fermionic. This $\Z_2$--grading should be also displayed by the symmetries of the physical system, that is, there are {\em even} transform
 ations which send bosons to bosons and fermions to fermions, and {\em odd} transformations, the so--called {\em supersymmetries}, which send bosons to fermions and vice versa. Mathematically speaking, implementing supersymmetry consists in passing from ordinary, ungraded algebraic structures such as vector spaces, Lie algebras etc, to their $\Z_2$--graded {\em super} counterparts, that is, to {\em super vector spaces}, {\em super Lie algebras}, etc. To see how this works concretely, we discuss supersymmetric $\sigma$--models next.
%
%
\subsection{Supersymmetric $\sigma$--models.}
Let $\Sigma$ be a Riemann surface (the ``string worldsheet''). A $\sigma$--{\em model} is a field theory whose bosonic fields consist of maps $\Sigma\to(M,g)$, where $M$ is a Riemannian manifold (the {\em target space} of the $\sigma$--model) --- for a precise and general definition of $\sigma$--models, see for instance the aforementioned~\cite{hal03}. For simplicity we assume $M$ to be compact, though this is not strictly necessary from a physical point of view. Depending on the amount of supersymmetry $\mc{N}$ one imposes (i.e.\ the labeled number of generators for the odd symmetry transformations), one speaks of an $\mc{N}=(1,1)$, $\mc{N}=(2,2)$ etc {\em supersymmetric} $\sigma$--model or target space, as $\mc{N}$ reflects the geometric structure of $M$. This is best illustrated by some examples.
\paragraph{$\mc{N}=(1,1)$ supersymmetry.} Let $g$ be a Riemannian metric, and consider the Hilbert space $\mc{H}$ given by completion of $\Omega^*(M)\otimes\C$ with respect to the canonical Hermitian structure induced by $g$. We have a natural $\Z_2$--grading on $\mc{H}$ coming from the decomposition into even and odd forms. The space of graded operators naturally carries a super Lie algebra structure whose bracket satisfies $[A,B]=(-1)^{\deg A\cdot\deg B+1}[B,A]$. The exterior differential $d:\Omega^{ev,od}(M)\otimes\C\to\Omega^{od,ev}(M)\otimes\C$ extends to a supersymmetry transformation. Further, the metric $g$ gives rise to a second supersymmetry transformation, namely the formal adjoint $d^*$. The resulting Laplacian $\Delta=dd^*+d^*\!d$ acts as a Hamiltonian which preserves the $\Z_2$--grading. Setting $Q_L=d+d^*$ and $Q_R=i(d-d^*)$ as in~\cite{wi82}, $Q_L$ and $Q_R$ satisfy satisfy $Q_L^2=Q_R^2=\Delta$. Further, $Q_L$ and $Q_R$ generate a vector space $\mf{g}^{odd}$ w
 hich we take as the odd part of a super vector space $\mf{g}$. The even part is spanned by $\Delta$. Since $d^2=d^{*2}=0$, the only non--trivial supercommutator relations are
$$
[Q_L,Q_L]=[Q_R,Q_R]=2\Delta.
$$
Hence, the supercommutator closes on $\mf{g}=\mf{g}^{ev}\oplus\mf{g}^{od}$ and therefore defines a subalgebra. In physicist's language, $\mf{g}$ is a $\mc{N}=(1,1)$--{\em supersymmetry algebra}\footnote{The notation $(1,1)$ refers to two fermionic fields on $\Sigma$, the so--called {\em left--} and {\em right--mover}, which makes that one can distinguish the supersymmetry generators $Q_L$ and $Q_R$ accordingly.}. 
\paragraph{$\mc{N}=(2,2)$ supersymmetry.} As we remarked above, on a K\"ahler manifold $(M^{2m},\omega,J)$ the complexified exterior algebra can be decomposed into $(p,q)$--forms $\Omega^{p,q}(M)$, and the derivative splits into $d=\partial+\overline{\partial}$. Taking formal adjoints, we end up with four supersymmetry transformations $G_1=\partial$, $G_1^*=\partial^*$, $G_2=i\overline{\partial}^*$ and $G_2^*=-i\overline{\partial}$. They satisfy
$$
[G_a,G_b]=[G^*_a,G^*_b]=0,\quad[G_a,G^*_b]=\delta_{ab}\Delta_{\partial},
$$
where $\Delta_{\partial}=\partial\partial^*+\partial^*\!\partial$ is the holomorphic Laplacian which coincides with the anti--holomor-phic Laplacian $\Delta_{\overline{\partial}}$ because of K\"ahlerness. In order to obtain a so--called $\mc{N}=(2,2)$ supersymmetry algebra, we have to enlarge this by the even {\em Lefschetz operator} $L=\wedge\omega:\Omega^*(M)\otimes\C\to\Omega^*(M)\otimes\C$ which is wedging with the K\"ahler form $\omega$, its formal adjoint $L^*$, and the {\em counting operator} $\Pi$ acting on $k$--forms by $(k-m){\rm Id}$. For these, the supercommutator relations are given by
$$
[\Pi,L]=2L,\quad[\Pi,L^*]=-2L^*,\quad[L,L^*]=\Pi
$$
that is, $L$, $L^*$ and $\Pi$ span an $\mf{sl}(2)$--subalgebra. Further, we have the {\em K\"ahler identities}
$$
\begin{array}{lcclcc}
\,[\Pi,G_1] & = & G_1,\quad&\,[\Pi,G_2] & = & -G_2,\\
\,[L,G_1] & = & 0,\quad&\,[L,G_2] & = & G_1,\\
\,[L^*,G_1] & = & G_2,\quad\,&[L^*,G_2] & = & 0.
\end{array}
$$
In particular, $G_1$ and $G_2$ transform as {\em spinors} under the $\mf{sl}(2)$--algebra action\footnote{Under the identification $SL(2,\C)\cong Spin(3,\C)$}. Mutatis mutandis, the same relations hold for $G_{1,2}^*$.

\begin{remark}
The operators of the $\mc{N}=(2,2)$ supersymmetry algebra can always be defined on a hermitian manifold, that is, a Riemannian manifold with compatible almost complex structure. However, the lack of integrability materialises in more complicated supercommutator relations. In particular, the super vector space $\mf{g}=\mf{g}^{ev}\oplus\mf{g}^{ev}$ is no longer closed under the supercommutator. On the other hand, K\"ahler manifolds are by no means the only geometric structures giving rise to $\mc{N}=(2,2)$--supersymmetry algebras. In~\cite{ghr84}, it was shown that the most general structure is a {\em bihermitian} manifold $(M,g,J_+,J_-,H)$, where $J_{\pm}$ are two almost complex structures preserving the metric $g$. Further, $H$ is a closed $3$--form which induces the twisted connections $\nabla^{\pm}=\nabla^g\pm H/2$. The integrability condition is
$$
\nabla^{\pm}J_{\pm}=0,\quad H=H^{1,2\pm}\oplus H^{1,2\pm}, 
$$
where the latter condition means that $H$ is the sum of a $(1,2)$ and a $(2,1)$--form with respect to {\em both} almost complex structures $J_+$ and $J_-$. A K\"ahler manifold arises as the special case $J_+=J=J_-$ and $H=0$. If $H=dB$, then we can interpret this special case as a $B$--{\em transform} of a K\"ahler structure, and more generally, bihermitian geometry as a {\em (twisted) generalised K\"ahler}--geometry (cf. Sections~\ref{btrafo} and~\ref{gengeo}).
\end{remark}

\noindent In view of the previous remark, it is enlightening to consider the case of a symplectic manifold.

\begin{example}
On a symplectic manifold $(M,\omega)$ we can always define the counting operator $\Pi$, the Lefschetz operator $L$ and its dual which we denote by $\Lambda$ in order to distinguish it from its formal adjoint $L^*$ (requiring a metric) considered earlier. Apart from the exterior differential $d$, we also have
$$
\widetilde{d}^*=[\Lambda,d]=\Lambda\circ d-d\circ\Lambda.
$$
This operator squares to zero and spans, together with $d$, the $2$--dimensional spin representation of the Lie algebra $\mf{sl}(2)$ generated by $\Pi$, $L$ and $\Lambda$, for
$$
\begin{array}{lcclcc}
\,[\Pi,d] & = & d,\quad&[\Pi,\widetilde{d}^*] & = & -\widetilde{d}^*\\
\,[L,d] & = & 0,\quad&[L,\widetilde{d}^*] & = & d,\\
\,[\Lambda,d] & = & \widetilde{d}^*,\quad&[\Lambda,\widetilde{d}^*] & = & 0.
\end{array}
$$
However, in absence of a Riemannian metric we do not dispose of any natural norm on $\Omega^*(M)\otimes\C$ in order to obtain a Hilbert space by completion. On the other hand, a classical result asserts that we can always choose a compatible almost--complex structure which by~(\ref{symcommet}) induces a Riemannian metric. In particular, we then have $\Lambda=L^*$. The metric also enables us to introduce a second $\mf{sl}(2)$--spin representation defined by the formal adjoints $d^*$ and $\widetilde{d}=(\widetilde{d}^*)^*$. But by lack of integrability, i.e.\ by the non--K\"ahlerness of the underlying almost K\"ahler structure, $\{d,\widetilde{d}\}$ fails to commute ($\{d,\widetilde{d}\}=0$ is in fact equivalent to K\"ahlerness), so that the super vector space spanned by these operators is not closed under the supercommutator. 
\end{example}
%
%
\subsection{Compactified type II supergravity}
For the moment being, we know five consistent supersymmetric string theories, namely {\em type I}, {\em heterotic} and {\em type II} string theory, where the latter two arise in two different flavours. Their low energy limit gives the corresponding supergravity theory, for instance type II supergravity on which we focus in this section. 

\bigskip

We first consider a ten dimensional space--time $M^{1,9}$, i.e.\ a Lorentzian manifold of dimension $10$. The bosonic fields come in two flavours; they are either NS-NS (NS=Neveu-Schwarz) or R-R (R=Ramond). To the former class belong
\begin{itemize}
\item the space--time metric $g$.
\item the dilaton field $\phi\in C^{\infty}(M)$.
\item the $H$--{\em flux} $H$, a closed $3$--form.
\end{itemize}
The R-R sector consists of a closed differential form $F$ of either even (type IIA) or odd degree (type IIB). The homogeneous components of $F$ are referred to as {\em Ramond--Ramond fields}. For sake of brevity, we shall only consider the case $F=0$, cf.~\cite{jewi05b} for the general case. Out of this datum one can build a Lagrangian $\mc{L}_{II}(g,\phi,H)$, and the critical points of the action functional $\mc{S}=\int_M\mc{L}_{II}$ determined by the Euler--Lagrange equations give the equations of motion\footnote{As usual in physics, one assumes all quantities to be decreasing fast enough at infinity so that the integral is well--defined.}. Further, one considers the action of a super Lie algebra whose action preserves the action functional and is generated by two {\em supersymmetry parameters} $\Psi_+$ and $\Psi_-$ (whence ``type II''), which are unit spinors of equal (type IIB) or opposite (type IIA) chirality. The model also contains fermionic fields, and these are suppo
 sed to be invariant under the supersymmetry transformations. This requirement is equivalent to the vanishing of the {\em supersymmetry variations} given by
\begin{eqnarray*}
\nabla_X\Psi_{\pm}\pm\frac{1}{4}(X\llcorner H)\cdot\Psi_{\pm} & = & 0\\
(d\phi\pm\frac{1}{2}H)\cdot\Psi_{\pm} & = & 0.
\end{eqnarray*}
The first pair of equations is referred to as the {\em gravitino equation} in the physics literature, while the latter is called {\em dilatino equation}. 

\bigskip

In order to find a solution, one usually makes a {\em compactification} ansatz, that is, one considers a space--time of the form $(M^{1,9},g^{1,9})=(\R^{1,3},g_0)\times(M^6,g)$, where $(\R^{1,3},g_0)$ is flat Minkowski space and $(M^6,g)$ a $6$--dimensional spinnable Riemannian manifold. This is not only a convenient mathematical ansatz, it also reflects the empiric fact that the phenomenologically tangible world is confined to three spatial dimensions plus time. Then, one tries to solve the above equations with fields living on $M^6$, trivially extended to the entire space--time. In essence, we are looking for a set $(g,\phi,H,\Psi_+,\Psi_-)$, that satisfies the gravitino and the dilatino equation on $M^6$. Simplifying further and setting $\phi=0$, $H=0$ and $\Psi_+=\Psi=\Psi_-$ (type IIB), we are left with the equation
$$
\nabla^g\Psi=0,
$$
which precisely characterises Calabi--Yau metrics on $M^6$. In the general case, the spinor fields $\Psi_+$ and $\Psi_-$ define each an underlying almost topological Calabi--structure $(M,\Omega_{\pm},\omega_{\pm})$ whose induced metrics coincide. This is reminiscent of the situation encountered in the discussion of the $\mc{N}=(2,2)$ supersymmetric $\sigma$--model characterised by a bihermitian geometry. In fact, we are again in territory of generalised geometry, and we will identify the geometry of type II supergravity compactifications as a generalised $SU(m)$--geometry. Again, a non--trivial $H$ reflects the lack of integrability of the underlying almost Calabi--Yau structure. More precisely,
$$
\nabla^g_X\Psi_{\pm}\pm\frac{1}{4}(X\llcorner H)\cdot\Psi_{\pm}=0
$$
boils down to $d\omega_{\pm}\not=0$, the right--hand side being determined by $H$, while we still have $d\Omega_{\pm}=0$ so that the underlying metric is still hermitian~\cite{chsa02}. 
%
%
%
%
%
\section{The B--field action}
\label{btrafo}
In light of the previous section, we wish to give geometric meaning to triples such as $(M,g,H)$ and $(M,\omega,H)$ etc. Here we enter the realm of {\em twisted generalised geometry} (cf. Section 7 in~\cite{hi03}). In this article, we consider a (mildly) simplified setup, where $H=dB$, which leads to the notion of a $B$--{\em field}, that is a $2$--form $B\in\Omega^2(M)$, and subsequently to the $B$--{\em transform} of a classical geometry.

\bigskip

The first problem which arises is that we cannot simply encode the additional datum of a $B$--field in a $G$--structure as in Section~\ref{classgeom}, as $B$ is an arbitrary $2$--form, and in particular, might have zeros. However, to induce a $G$--structure requires the tensors to enjoy special global algebraic properties. Rather, one associates a transformation with $B$ which acts as the identity on the zero locus of $B$. Concretely, consider $B$ as a map $TM\to T^*\!M$. Then extend this to an endomorphism of $\E$ which we write as a block matrix
$$
B\leadsto\left(\begin{array}{cc} 0 & 0\\B & 0\end{array}\right).
$$
Taking the usual matrix exponential yields the {\em isomorphism}
$$
e^B=\left(\begin{array}{cc} Id & 0\\B & Id\end{array}\right),
$$
that is, $X\oplus\xi\in\E$ is sent to $X\oplus\big(B(X)+\xi\big)\in\E$.
In this way, $B$ does not only induce an isomorphism, but in fact an orientation--preserving isometry of $\E$. Indeed, it preserves the natural orientation on $\E$ and the inner product of signature $(n,n)$ given by
$$
(X\oplus\xi,X\oplus\xi)=\xi(X).
$$
In particular, $\E$ is associated with an $SO(n,n)$--structure which can be further reduced to $SO(n,n)_0$ if, as we always assume, $M$ is orientable\footnote{Here and in the sequel, $G_0$ denotes the identity component of a Lie group $G$.}.

\bigskip

Further, $B$ also acts on forms of mixed degree $\rho\in\Omega^{ev,od}(M)$, namely by
$$
e^B\wedge\rho=(1+B+\frac{1}{2}B\wedge B+\ldots)\wedge\rho\in\Omega^{ev,od}(M).
$$
If $B$ is a {\em closed} $B$--field, then it commutes with the natural differential operators defined on $\Omega^*(M)$ and $\E$, namely the exterior differential $d:\Omega^{ev,od}(M)\to\Omega^{od,ev}(M)$ and the {\em Courant bracket}~\cite{co90}
$$
\llbracket X\oplus\xi,Y\oplus\eta\rrbracket=[X,Y]+\mc{L}_X\eta-\mc{L}_Y\xi-\frac{1}{2}d(X\llcorner\eta-Y\llcorner\xi)
$$
($\mc{L}$ denoting Lie derivative). This bracket is skew--symmetric, but does not satisfy the Jacobi identity. 

\bigskip

The interpretation of a ``classical geometry'' with a $B$--field consists in thinking of the resulting structure as being obtained by a $B$--field action on the original setup, much in the vein of the remark at the end of Section~\ref{classgeom}, where we considered the $f$--transform of a(n) (almost) symplectic manifold. In this picture, {\em closed} $B$--fields correspond to diffeomorphisms. Consequently, we need to rephrase classical geometries in terms of even or odd forms or additional algebraic datum on the bundle $\E$, together with an integrability condition involving $d$ or $\llbracket\cdot\,,\cdot\rrbracket$. Defining the resulting $B$--transforms in an invariant way gives rise to so--called {\em generalised geometries}, where we have effectively enhanced the transformation group to $\diff(M)\ltimes\Omega^2_{cl}(M)$. We discuss examples next.
%
%
%
%
%
\section{Generalised geometry}
\label{gengeo}
\paragraph{Generalised Riemannian metrics~\cite{gu04},~\cite{wi05}.}
Let $(M,g)$ be a Riemannian manifold.\hspace{-2pt} Viewing $g$ as an isomorphism $T\!M\!\!\to\!T^*\!M$, it is lifted to $T\!M\oplus T^*\!M$ by
$$
\mc{G}_g=\left(\begin{array}{cc} 0 & g^{-1}\\g & 0\end{array}\right).
$$
A $B$--field acts via
$$
\mc{G}_{g,B}=e^B\mc{G}_g e^{-B}=\left(\begin{array}{cc} -g^{-1}B & g^{-1}\\ g-Bg^{-1}B & Bg^{-1}\end{array}\right).
$$
So, symbolically, a {\em Riemannian manifold with $B$--field} is pictured by
$$
\shabox{$(M,g)\longleftrightarrow(M,\mc{G}_g)\quad\stackrel{e^B}{\leadsto}\quad(M,g,B)\longleftrightarrow(M,\mc{G}_{g,B}=e^B\mc{G}_ge^{-B})$.}
$$
As invariant definition we adopt the following.

\begin{definition}
An endomorphism $\mc{G}$ of $\E$ is said to be a {\em generalised Riemannian metric} if $\mc{G}^2=Id$ and $\big(\mc{G}(Z),Z\big)>0$ for all $0\not=Z\in\E$. 
\end{definition}

\noindent The endomorphisms $\mc{G}_g$ and $\mc{G}_{g,B}$ are special cases of this. To see what a generalised Riemannian metric means in terms of structure groups, we remark that defining $\mc{G}$ is equivalent to decomposing $\E=V^+\oplus V^-$ into a maximally positive and negative definite subbundle, corresponding to the $\pm1$--eigenbundles of $\mc{G}$. Let $g_{\pm}=(\cdot\,,\cdot)_{V^{\pm}}$ denote the resulting bundle metric on $V^{\pm}$. At a point $x\in M$, this splitting is preserved under the action of $SO(n,n)$ by the subgroup $SO(V^+_x,g_{+x})\times SO(V^-_x,g_{-x})$. Locally, this yields the family of representations
$$
\begin{array}{ll}
r_x:(A,B)\in SO(n)\times SO(n)\mapsto\left(\begin{array}{cc}A&0\\0&B\end{array}\right)\!\!\!\!& \in  SO(V^+_x,g_{+x})\times SO(V^-_x,g_{-x}) \\
\!\!\!\!&\subset SO(T_xM\oplus T^*_x\!M)_0,
\end{array}
$$
where the block matrix is written with respect to the decomposition $V^+_x\oplus V^-_x$. Hence, a generalised Riemannian metric is associated with an $SO(n)\times SO(n)$--structure. Since the isotropic subbundles $TM$ and $T^*\!M$ intersect the positive and negative subbundles $V^+$ and $V^-$ trivially, the projections $\pi_{\pm}:TM\subset\E\to V^{\pm}$ are isomorphisms. Further, one can show that there exists a unique Riemannian metric $g$ and a $2$--form $B$ such that \begin{equation}\label{isom}
\pi_{\pm}(X)=X\oplus(\pm g+B)(X).
\end{equation} 
Then, the generalised Riemannian metric corresponding to $V^+\oplus V^-$ is just $\mc{G}_{g,B}$, and $\pi_{\pm}:(T,\pm g)\to (V^{\pm},g_{\pm})$ are promoted to isometries. In particular, any $\mc{G}$ is of the form $\mc{G}_{g,B}$, so that modulo arbitrary $B$--field transformations, we have only one generalised Riemannian metric for every usual Riemannian metric. However, transforming with a {\em non--closed} $B$--field modifies the Levi--Civita connection $\nabla^g$ which gets replaced by the two connections $\nabla^{\pm}=\nabla^g\pm dB/2$ (cf. also Theorems~\ref{biherm} and~\ref{susyeq}). 
\paragraph{Generalised complex structures~\cite{gu04},~\cite{hi03}.}
Next consider a complex manifold $(M,J)$. We lift this to $\E$ by
\begin{equation}\label{gencomJ}
\mc{J}_J=\left(\begin{array}{cc} -J & 0\\0 & J^*\end{array}\right).
\end{equation}
Its $B$--field transform is
$$
\mc{J}_{J,B}=e^B\mc{J}_Je^{-B}=\left(\begin{array}{cc}-J & 0\\-BJ-J^*B & J^*\end{array}\right),
$$
so a {\em almost complex manifold with $B$--field} is described by
$$
\shabox{$(M,J)\longleftrightarrow(M,\mc{J}_J)\quad\stackrel{e^B}{\leadsto}\quad(M,J,B)\longleftrightarrow(M,\mc{J}_{J,B}=e^B\mc{J}_Je^{-B})$.}
$$
The invariant definition is this.

\begin{definition}
An {\em almost generalised complex structure} $(M,\mc{J})$ is given by $\mc{J}\!\in\End(\E)$ such that $\mc{J}^2=-Id$ and $\mc{J}$ is an isometry for $(\cdot\,,\cdot)$. 
\end{definition}

The stabiliser group of an almost complex structure compatible with an inner product of split signature inside $SO(2m,2m)$ is $U(m,m)$, so $(M,\mc{J})$ is associated with an $U(m,m)$--structure. As in the classical case, an almost generalised complex manifold $M$ must be even--dimensional. This notion comprises $\mc{J}_J$ and $\mc{J}_{J,B}$ as special instances, but is far richer.

\begin{examples}\hfill\newline
\indent(i) {\em Symplectic manifolds.} An example which is not of the form $\mc{J}_{J,B}$ is provided by a symplectic manifold $(M,\omega)$. It induces an almost generalised complex structure by
\begin{equation}\label{gencomsymp}
\mc{J}_{\omega}=\left(\begin{array}{cc}0&\omega^{-1}\\\omega&0\end{array}\right).
\end{equation}
\indent(ii) {\em Products}. Two almost generalised complex manifolds $(M_{1,2},\mc{J}_{1,2})$ induce a natural almost generalised complex structure $\mc{J}_1\oplus\mc{J}_2$ on the product $T(M_1\times M_2)\cong TM_1\oplus TM_2$.
\end{examples}

\noindent An almost generalised complex structure is also equivalent to a decomposition $(\E)\otimes\C=W_{\mc{J}}\oplus\overline{W_{\mc{J}}}$ into isotropic complex subbundles. For example, we find $W_J=T^{0,1}\oplus T^{1,0*}$ for $\mc{J}_J$ and $W_{\omega}=\{X\oplus i\omega(X,\cdot)\}\,|\,X\in\Gamma(TM\otimes\C)\}$ for $\mc{J}_{\omega}$.

\begin{definition}
We speak of an {\em integrable} or simply a {\em generalised complex structure} $(M,\mc{J})$, if the generalised Nijenhuis tensor
$$
\mc{N}^{\mc{J}}(Z_0,Z_1)=\llbracket Z_0,Z_1\rrbracket -\llbracket \mc{J}Z_0,\mc{J}Z_1\rrbracket +\mc{J}(\llbracket Z_0,\mc{J}Z_1\rrbracket +\llbracket \mc{J}Z_0,Z_1\rrbracket ) 
$$
vanishes for all $Z_0,\,Z_1\in\Gamma(\E)$.
\end{definition}

Alternatively, we can require $W_{\mc{J}}$ to be closed under the Courant bracket. In particular, an integrable structure naturally transforms under closed $B$--fields. In the case of $\mc{J}_J$ and $\mc{J}_{\omega}$, integrability is equivalent to the integrability of the underlying classical structure, i.e.\ $(M,J)$ is a complex manifold and $(M,\omega)$ is symplectic. The corresponding form picture is less obvious, as forms of homogeneous degree are not acted on by $B$--fields --- here is where $\Omega^{ev,od}(M)$ comes in. We first define an action $\bullet:(\E)\times\Omega^{ev,od}(M)\to\Omega^{od,ev}(M)$ by
\begin{equation}\label{clifford}
(X\oplus\xi)\bullet\rho=-X\llcorner\rho+\xi\wedge\rho.
\end{equation}
This squares to $-(X\oplus\xi,X\oplus\xi){\rm Id}$ and naturally extends to $\Lambda^*(\E)$. In particular, the annihilator $W_{\rho}=\{Z\in(\E)\otimes\C\,|\,Z\bullet\rho=0\}$ of an even or odd complex form $\rho\in\Omega^{ev,od}(M)\otimes\C$ is totally isotropic. We call $\rho$ {\em pure} if $W_{\rho}$ is a subbundle of maximal rank, in analogy with pure spinor fields. A pure form is necessarily even or odd. Further, there is a 1-1--correspondence between line bundles in $\Lambda^*T^*\!M\otimes\C$ whose local trivialisations consist of pure spinors, and maximally isotropic subbundles of $(\E)\otimes\C$. The line bundle corresponding to $W_{\mc{J}}$ is referred to as the {\em canonical bundle} of $\mc{J}$ (the terminology will be justified below) and is written $\mc{K}_{\mc{J}}$. This yields a decomposition
$$
\Omega^*(M)=\oplus_{k=0}^n U_k,\mbox{ where } U_0=\mc{K}_{\mc{J}},\,U_k=\Lambda^kW^*_{\mc{J}}\bullet U_0.
$$
As $d$ changes the parity of a form, we can define $\partial=\pi_{k-1}\circ d:\Gamma(U_k)\to\Gamma(U_{k-1})$ and $\overline{\partial}=\pi_{k+1}\circ d:\Gamma(U_k)\to\Gamma(U_{k+1})$. Then $\mc{J}$ is integrable if and only if $d=\partial+\overline{\partial}$.
Further, there are also special coordinates for generalised complex structures, although this requires some care. We define $t(x)$, the {\em type of $\mc{J}$ at $x$}, to be the dimension of $\pi(W_{\mc{J},x})\subset T_xM\otimes\C$, the projection of $W_x$ to $T_xM$. If it is locally constant around $x$ (which does not always hold), then some neighbourhood of $x$ is equivalent, modulo diffeomorphism and closed $B$--field transformation, to an open neighbourhood $U_1\times U_2$ in $\C^t\times\R^{2(m-t)}$ with induced generalised complex structure coming from the standard $t$--dimensional complex manifold $\C^t$ times the standard $2(m-t)$--dimensional real symplectic manifold $(\R^{2(m-t)},\omega_0)$.
\paragraph{Generalised K\"ahler structures~\cite{gu04}.} Next we generalise the notion of a classical K\"ahler manifold.

\begin{definition}
A {\em generalised K\"ahler structure} on $M$ is defined by a pair of commuting generalised complex structures $(\mc{J}_0,\mc{J}_1)$ such that $\mc{G}=-\mc{J}_0\mc{J}_1$ defines a generalised Riemannian structure.
\end{definition}

\begin{examples}\hfill\newline
\indent(i) Consider a classical K\"ahler manifold $(M,\omega,J)$. The generalised complex structures $\mc{J}_J$ and $\mc{J}_{\omega}$ in~(\ref{gencomJ}) and~(\ref{gencomsymp}) commute indeed. Further,
$$
-\mc{J}_J\mc{J}_{\omega}=\left(\begin{array}{cc}0&g^{-1}\\g&0\end{array}\right)=\mc{G}_0,
$$
where $g:TM\to T^*\!M$ is the K\"ahler metric given by~(\ref{symcommet}). Therefore, a {\em K\"ahler manifold with $B$--field} is
$$
\shabox{$(M,J,\omega)\longleftrightarrow(M,\mc{J}_J,\mc{J}_{\omega})\quad\stackrel{e^B}{\leadsto}\quad(M,J,\omega,B)\longleftrightarrow(M,\mc{J}_{J,B},\mc{J}_{\omega,B}$).}
$$

(ii) In~\cite{hi06}, compact examples which are not a closed $B$--transform of a K\"ahler structure as in (i) were constructed on $\C P^2$ and $\C P^1\times\C P^1$.
\end{examples}

In terms of bundles, a generalised K\"ahler manifold can be characterised as follows. Since $\mc{J}_0$ and $\mc{J}_1$ commute, the decomposition $(\E)\otimes{\C}=W_0\oplus\overline{W_0}$,
where $W_0=W_{\mc{J}_0}$, is stable under $\mc{J}_1$. This implies a further decomposition of $W_0$ into the $\pm i$--eigenspaces $W^{\pm}_0$ of $\mc{J}_1$. Hence, $(\E)\otimes\C=W^+_0\oplus W^-_0\oplus\overline{W^+_0}\oplus\overline{W^-_0}$.
Moreover, since $\mc{G}=-\mc{J}_0\mc{J}_1$,
\begin{equation}\label{genkahl}
V^{\pm}\otimes\C=W^{\pm}_0\oplus\overline{W^{\pm}_0}.
\end{equation}
Indeed, $\mc{J}_0\mc{J}_1w=\mp w$ for $w\in W^{\pm}_0$ etc. One can then prove

\begin{proposition}\label{bundlecrit}
An almost generalised K\"ahler structure $(M^{2m},\mc{J}_0,\mc{J}_1)$ is equivalent to specifying two complex rank $m$ subbundles $W^{\pm}_0$ of $(\E)\otimes\C$ satisfying
\begin{itemize}
\item $W_0^{\pm}$ is isotropic
\item $W^+_0\perp W^-_0,\overline{W^-_0}$
\item $W^{\pm}_0$ is definite in the sense that $\pm(Z,\overline{Z})>0$ for all $Z\in W^{\pm}_0$ unless $Z=0$.
\end{itemize}
In this picture, integrability is equivalent to each of the subbundles $W_0^+$, $W_0^-$ and $W_0^+\oplus W_0^-$ being closed under the Courant bracket.
\end{proposition}

The structure group point of view gives yet another way of analysing generalised K\"ahler structures. The decomposition~(\ref{genkahl}) defines almost complex structures $\mc{J}_{\pm}$ on $V^+$ and $V^-$. Further, this structure is compatible with the positive and negative metric $g_{\pm}$ on $V^{\pm}$. Hence $V^+$ and $V^-$ are associated with a $U(m)$--structure, and the structure group is therefore $U(m)\times U(m)$. Further, these K\"ahler structures can be pulled back to $TM$ via the isomorphisms $\pi_{\pm}$ of~(\ref{isom}), giving rise to {\em two} K\"ahler structures $(\omega_{\pm},J_{\pm})$ on $M$ which share the same K\"ahler metric. Therefore, these K\"ahler structures can be characterised by $(M,g,J_{\pm})$.

\begin{theorem}\label{biherm}
An almost generalised K\"ahler structure $(\mc{J}_0,\mc{J}_1)$ corresponds to the datum $(g,B,J_+,J_-)$, where $J_{\pm}$ are two (and in general not commuting) almost complex structures on $M$, $g$ a Riemannian metric which is hermitian with respect to both $J_+$ and $J_-$, and $B\in\Omega^2(M)$, coming from the inclusions
$$
\begin{array}{ccccc}
U(m)\times U(m) & \subset & SO(2m)\times SO(2m) & \subset & SO(2m,2m)_0.\\[5pt]
(J_+,J_-) & & (g,B) & &
\end{array}
$$
In this picture, the almost K\"ahler structure is integrable, if and only if
$J_{\pm}$ are integrable complex structures, and the corresponding K\"ahler forms satisfy
$$
d^c_-\omega_-=dB=-d^c_+\omega_+,
$$
where $d^c_{\pm}=i(\overline{\partial}_{\pm}-\partial_{\pm})$. This can be rephrased in terms of connections by saying that
$$
\nabla^{\pm}J^{\pm}=0,\quad dB\in\Omega^{1,2}(M)\oplus\Omega^{2,1}(M),
$$
where $\nabla^{\pm}=\nabla^g\pm dB/2$.
\end{theorem}

\begin{remark}
In the more general case of a {\em twisted} generalised K\"ahler structure, $dB$ gets replaced by an arbitrary closed $3$--form $H$ in Theorem~\ref{biherm}. This gives precisely the conditions found on an $\mc{N}=(2,2)$--supersymmetric target space described in Section~\ref{physicsI}.
\end{remark}
\paragraph{Generalised Calabi--Yau structures~\cite{hi03}.}
A {\em generalised Calabi--Yau structure}~\cite{hi03} mimics the classical relation between complex and special complex manifolds in the generalised setting, rather than generalising classical Calabi--Yau manifolds (which effectively leads to the notion of a generalised $SU(m)$--structure, cf. next paragraph)\footnote{Using our jargon consistently, we would rather refer to this structure as a {\em generalised special complex structure}, but we prefer to stick to Hitchin's original terminology to avoid all too certain confusion.}. As on the classical level, this involves differential forms, of mixed even or odd degree, with special algebraic properties. 

\bigskip

Let us introduce on $\Omega^*(M)$ the $\Omega^{2m}(M)$--valued bilinear form 
\begin{equation}\label{innprod}
\langle\rho,\tau\rangle=[\rho\wedge\widehat{\tau}]^{2m}.
\end{equation}
Here, $[\,\cdot\,]^{2m}$ denotes projection on the top degree and $\hat{\cdot}$ acts on homogeneous components of degree $p$ by multiplication with $(-1)^{p(p+1)/2}$. It is symmetric for $m$ even and skew for $m$ odd.

\begin{definition}
An {\em almost generalised Calabi--Yau structure} on $M^{2m}$ is defined by an even or odd pure form $\rho\in\Omega^{ev,od}(M)\otimes\C$  with $\langle\rho,\overline{\rho}\rangle\not=0$. The structure $(M,\rho)$ is referred to as {\em integrable} or simply as {\em generalised Calabi--Yau} if $d\rho=0$.
\end{definition}

\noindent Purity implies that the annihilator $W_{\rho}$ is of maximal dimension. The second condition is equivalent to $W_{\rho}\cap W_{\overline{\rho}}=\{0\}$~(\cite{ch96}, p.143). Since $W_{\overline{\rho}}=\overline{W_{\rho}}$, an almost generalised Calabi--Yau structure induces a decomposition $(\E)\otimes\C=W_{\rho}\oplus\overline{W_{\rho}}$, hence an almost generalised complex structure. It is integrable if $d\rho=0$. Put differently, an almost generalised Calabi--Yau manifold is an almost generalised complex manifold whose canonical line bundle $\mc{K}_{\mc{J}}$ is trivialised by a global section $\rho$ which is ``generalised holomorphic'' in the sense that $\overline{\partial}\rho=0$.

\begin{examples}\hfill\newline\indent
(i) A symplectic manifold gives rise to the even closed form $\rho=e^{i\omega}$. Then $\langle\rho,\overline{\rho}\rangle=c\omega^m$ for some non--zero constant. Furthermore, $W_{\rho}=W_{\omega}$, so $\rho$ is pure. The induced generalised complex structure is $\mc{J}_{\omega}$. A {\em symplectic manifold with $B$--field} is thus
$$
\shabox{$(M,\omega)\longleftrightarrow(M,e^{i\omega})\quad\stackrel{e^B}{\leadsto}\quad(M,\omega,B)\longleftrightarrow(M,e^{i\omega+B})$.}
$$

(ii) A special complex manifold $(M,\Omega)$ induces a generalised Calabi--Yau structure by $\rho=\Omega$. Here, $\langle\Omega,\overline{\Omega}\rangle=\pm\Omega\wedge\overline{\Omega}\not=0$ by the very definition of $\Omega$. Furthermore, if $J$ denotes the induced complex structure, $W_{\rho}=T^{0,1}\oplus T^{1,0*}$. In particular, the induced generalised complex structure is $\mc{J}_J$. A {\em special complex manifold with $B$--field} is then simply
$$
\shabox{$(M,\Omega)\quad\stackrel{e^B}{\leadsto}\quad(M,\Omega,B)\longleftrightarrow(M,e^B\wedge\Omega)$.}
$$

(iii) A classical Calabi--Yau manifold $(M,\omega,\Omega)$ gives rise to the generalised Calabi--Yau structures defined by $\exp(i\omega)$ and $\Omega$.

(iv) The $B$--field can be thought of as interpolating between these two extreme cases: Consider the case of a holomorphic symplectic manifold $(M^{2k},\omega_c=\omega_1+i\omega_2)$ of complex dimension $\dim_{\C}M=2k$. Then $\omega_{1,2}$ are real symplectic forms and $\Omega^{2k}=\omega^k_c/k!$ makes $M^{4k}$ (now seen as a real manifold) special complex. Consider for $t\not=0$ the family of generalised Calabi--Yau structures 
$$
\rho_t=t^k\cdot e^{\omega_c/t}=e^{\omega_1/t}\wedge(t^k\cdot e^{i\omega_2/t})=t^k+\ldots+\frac{1}{k!}\omega^k_c.
$$
Then $\rho_t$ is the $B$--field transform of $\big(M^{4k},t^k\cdot\exp(i\omega_2/t)\big)$ which converges to $(M^{4k},\Omega^{2k})$ for $t\to0$. 
\end{examples}

In terms of structure groups, a pure form $\rho$ trivialises the canonical line bundle $\mc{K}_{\mc{J}}$ of the almost generalised complex structure it induces. It remains to see that it also trivialises $\Lambda^{2m}W^*_{\rho}$ to conclude that an almost generalised Calabi--Yau manifold has structure group $SU(m,m)$. For this, we need to to make contact with $\E$--spinor fields. Here, 
$$
\mb{S}(\E)_{\pm}=\big(\Lambda^{ev,od}T^*\!M\otimes\sqrt{\Lambda^nTM}\big)\otimes\C,
$$
where Clifford multiplication is given by $\bullet$ as in~(\ref{clifford}) acting on the form part. As $M$ is orientable, we can choose a trivialisation of $\Lambda^nTM$ and identify $\E$--spinor fields with complex differential forms. Any scale--invariant property such as purity is not affected by this choice and makes therefore sense on both sides. As for $TM$--spinor fields, there is an isomorphism
\begin{equation}\label{genfierzing}
\mb{S}(\E)\otimes\mb{S}(\E)\cong\Lambda^*(\E)\otimes\C
\end{equation}
which by III.3.2~\cite{ch96} implies $(\mc{K}_{\mc{J}}\otimes\sqrt{\Lambda^{2m}TM})^2\cong\Lambda^{2m}W_{\rho}$. As $\mc{K}_{\mc{J}}$ is trivial, so is $\Lambda^{2m}W_{\rho}$ and thus $\Lambda^{2m}W_{\rho}^*$. Consequently, the structure group is $SU(m,m)$, which can be also thought of as the stabiliser of $\rho_x$ under the action of the spin group $Spin(2m,2m)$ on $\mb{S}(\E)_{\pm\,x}$.
\paragraph{Generalised $SU(m)$--structures~\cite{gmwi06},~\cite{jewi05},~\cite{jewi05b},~\cite{wi06b}.}
Next we wish to implement the counterpart of a classical Calabi--Yau manifold in the generalised setting.

\begin{definition}
An {\em almost generalised $SU(m)$--structure} $(M,\rho_0,\rho_1)$ is given by a pair of forms such that
\begin{itemize}
	\item $\rho_0$ and $\rho_1$ induce each a generalised Calabi--Yau structure such that $\langle\rho_0,\overline{\rho_0}\rangle=c\langle\rho_1,\overline{\rho_1}\rangle$ for some constant $c$.
	\item the induced generalised complex structures $(M,\mc{J}_{\rho_0},\mc{J}_{\rho_1})$ define an almost generalised K\"ahler structure.
\end{itemize}
We speak of an {\em integrable} or simply of a {\em generalised $SU(m)$--structure} if both $\rho_0$ and $\rho_1$ are closed\footnote{In~\cite{gu04}, this type of structure was referred to as a {\em generalised Calabi--Yau metric}. In the light of Theorem~\ref{susyeq}, we prefer the term generalised $SU(m)$--structure which is in line with the notion of a generalised $G_2$--structure~\cite{wi06}.}.
\end{definition}

\begin{example}
Take a classical Calabi--Yau manifold $(M,\omega,\Omega)$ and let $\rho_0=\exp(i\omega)$, $\rho_1=\Omega$. From the definition of a Calabi--Yau, $\Omega\wedge\overline{\Omega}=(-1)^{m(m-1)/2}i^m\omega^m$. Further, as we have seen above, $(\mc{J}_{\rho_0}=\mc{J}_{\omega},\mc{J}_{\rho_1}=\mc{J}_{\Omega})$ is generalised K\"ahler. Hence, any classical Calabi--Yau manifold induces a canonical generalised $SU(m)$--structure (but two different generalised Calabi--Yau structures). A {\em Calabi--Yau manifold with $B$--field} is thus
$$
\shabox{$(M,\Omega,\omega)\longleftrightarrow(M,\Omega,e^{i\omega})\quad\stackrel{e^B}{\leadsto}\quad(M,\Omega,\omega,B)\longleftrightarrow(M,e^B\wedge\Omega,e^{i\omega+B})$.}
$$
\end{example}

Whether or not $(\rho_0,\rho_1)$ induces a generalised K\"ahler structure can be read off directly from the forms. If $W_0$ and $W_1$ are the annihilators of $\rho_0$ and $\rho_1$ in $(\E)\otimes\C$, then $W_0^+=W_0\cap W_1$ and $W_0^-=W_0\cap\overline{W_1}$ satisfy the first two properties of Proposition~\ref{bundlecrit}. A necessary and sufficient criterion for $W_0^{\pm}$ to be of complex rank $m$ is given by Proposition III.4.4 in~\cite{ch96}: Regarding $\rho_0$ and $\rho_1$ as spinors and using the identification~(\ref{genfierzing}), this is equivalent for the forms $\rho_1\otimes\rho_2$ and $\rho_1\otimes\overline{\rho_2}$ to have no homogeneous components of degree strictly less than $m$ (this condition is clearly scale--invariant). 
From this we deduce the

\begin{proposition}
A pair $(\rho_0,\rho_1)$ of pure forms defines a generalised $SU(m)$--structure if and only if the following conditions hold:
\begin{itemize}
\item $\langle\rho_0,\overline{\rho}_0\rangle=c\langle\rho_1,\overline{\rho}_1\rangle>0$ for some constant $c$.
\item the forms $\rho_0\otimes\rho_1,\,\rho_0\otimes\overline{\rho_1}\in\Lambda^*(\E)\otimes\C$ have no homogeneous components of degree strictly less than $m$.
\item for all $Z$ in $W_0^+=W_0\cap W_1$ or $W_0^-=W_0\cap\overline{W_1}$, we have $(\overline{Z},Z)<0$.
\end{itemize} 
\end{proposition}

To derive the structure group point of view, we start with the decomposition $V^{\pm}\otimes\C=W_0^{\pm}\oplus\overline{W_0^{\pm}}$. The canonical line bundles of the $U(m)$--structures on $V^{\pm}$ are therefore $\kappa_{\pm}=\Lambda^mW^{\pm *}_0$.
On the other hand, 
$$
\mc{K}_0=\Lambda^{2m}W_0^*=\Lambda^{2m}(W^+_0\oplus W^-_0)^*\cong\Lambda^mW^{+*}_0\otimes\Lambda^mW^{-*}_0=\kappa_+\otimes\kappa_-.
$$
Since $\mc{J}_1=\mc{G}\mc{J}_0$ acts on $W_0^+\oplus\overline{W^-_0}$ as $i\cdot Id$, we deduce that
$$
\mc{K}_1=\Lambda^{2m}W_1^*=\Lambda^{2m}(W^+_0\oplus \overline{W_0^-})^*\cong\kappa_+\otimes\overline{\kappa_-}.
$$
Since $\rho_0$ and $\rho_1$ trivialise $\mc{K}_0$ and $\mc{K}_1$ respectively, $\kappa_+\cong\overline{\kappa_-}\cong\kappa_-$, hence $\kappa_+$ and $\kappa_-$ are trivial.
Consequently, $V^+$ and $V^-$ are associated with $SU(m)$--structures, and the structure group is therefore $SU(m)\times SU(m)$. Following the generalised K\"ahler case, we pull these back to the tangent bundle and obtain two $SU(m)$--structures on $TM$ given by $(\Omega_{\pm},\omega_{\pm})$. These induce the same metric and can be thus equivalently characterised by $(g,\Psi_{\pm})$ with two pure $TM$--spinor fields $\Psi_+$ and $\Psi_-$ (cf. Section~\ref{classgeom}). Conversely, given $(g,B,\Psi_+,\Psi_-)$, we can define $\rho_0=e^B\wedge\mc{A}(\Psi_+)\otimes\Psi_-$ and $\rho_1=e^B\wedge\Psi_+\otimes\Psi_-$. One can show that $(\rho_0,\rho_1)$ defines indeed an almost generalised Calabi--Yau structure. Further, so does $(e^{-\phi}\cdot\rho_0,e^{-\phi}\cdot\rho_1)$ for any function $\phi\in C^{\infty}(M)$, as purity is not affected by rescaling.

\begin{theorem}\label{susyeq}
For any almost generalised $SU(m)$--structure given by $(M,\rho_0,\rho_1)$, there exists a metric $g$, a $2$--form $B$, a function $\phi$ and two pure $TM$--spinor fields $\Psi_+$ and $\Psi_-$ such that
$$
\rho_0=e^{-\phi}\cdot e^B\wedge\mc{A}(\Psi_+)\otimes\Psi_-,\quad\rho_1=e^{-\phi}\cdot e^B\wedge\Psi_+\otimes\Psi_-.
$$
This datum corresponds to the inclusions
$$
\begin{array}{ccccc}
SU(m)\times SU(m) & \subset & Spin(2m)\times Spin(2m) & \subset & Spin(2m,2m).\\[5pt]
(\Psi_+,\Psi_-) & & (g,B) & & 
\end{array}
$$
and
$$
\begin{array}{ccc}
Spin(2m,2m) & \subset &  \R_{>0}\times Spin(2m,2m).\\[5pt]
\phi & &
\end{array}
$$
Further, the structure is integrable, that is, $d\rho_0=0$ and $d\rho_1=0$, if and only if
$$
\nabla^g_X\Psi_{\pm}\pm\frac{1}{4}(X\llcorner dB)\cdot\Psi_{\pm}=0,\quad(d\phi\pm\frac{1}{2}dB)\cdot\Psi_{\pm}=0.
$$
\end{theorem}

\begin{remark}
Again, there is a twisted notion of a generalised $SU(m)$--structure which gives the general case for closed, but not exact $H$, for which Theorem~\ref{susyeq} holds with $H$ in place of $dB$. As a result, we have identified the geometric structure underlying the type II supergravity compatifications discussed in Section~\ref{physicsI} as generalised $SU(3)$--structures. 
\end{remark}

\begin{example}
For $\Psi_+=\Psi_-=\Psi$, $B=0$ and $\phi=0$, we recover the usual Calabi--Yau manifold as given by~(\ref{CYspinor}).
\end{example}

\noindent This formulation yields easily the Ricci and scalar curvature $\ric^g$ and $S^g$ of a generalised $SU(m)$--metric $g$~\cite{jewi05b}.

\begin{corollary}
Let $g$ be the metric of an integrable generalised $SU(m)$--structure. Then
$$
\ric^g(X,Y)=-2H^{\phi}(X,Y)+\frac{1}{2}g(X\llcorner dB,
Y\llcorner dB)
$$
and
$$
S^g=2\Delta\phi+\frac{3}{2}|dB|^2,
$$
where $H^{\phi}(X,Y)=X.Y.\phi-\nabla^g_XY.\phi$ is the Hessian of $\phi$.
\end{corollary}

\noindent A further consequence are two striking {\em no--go theorems}, reflecting the fact that the equations $d\rho_0=0$ and $d\rho_1=0$ are heavily overdetermined.

\begin{corollary}\hfill\newline
\indent(i) If $M$ is compact, then $dB=0$.

(ii) If $\phi\equiv const$, then $dB=0$.

In either case, a generalised $SU(m)$--structure is equivalent to a usual Calabi--Yau manifold via a closed $B$--field transformation.
\end{corollary}

\begin{remark}
This is in line with similar statements for {\em generalised $G_2$-- and $Spin(7)$--structures} \cite{wi06} which admit $TM$--spinorial characterisations in the vein of Theorem~\ref{susyeq}, but stands in sharp contrast to generalised K\"ahler structures where compact examples do exists, as pointed out above. Note that these no--go theorems only hold for the closed case $d\rho=0$, $d\rho_1=0$. For the non--closed case, compact non--classical examples with $dB\not=0$ might well exist (cf. for instance~\cite{fito06} for a compact generalised $G_2$--example). This kind of inhomogeneous conditions arise naturally when considering constrained variational problems or type II compactifications with R-R fields (cf. Section~\ref{physicsI}).
\end{remark}
%
%
%
%
%

\end{document}